\documentstyle[12pt]{article}
\begin{document}

%LATEX 209

%ALPHABET...
\def\c{{\bf C}}
\def\n{{\bf N}}
\def\z{{\bf Z}}
\def\u{{\bf U}}
\def\f{{\bf F}}
\def\t{{\bf T}}
\def\s{{\bf S}}
\def\o{{\bf O}}
\def\w{{\bf W}}
\def\g{{\bf G}}
\def\r{{\bf R}}
\def\p{{\bf P}}
\def\pp{{\bf p}}
\def\ll{{\bf L}}
\def\l{{\cal L}}
\def\co{{\cal C}}

%CONJUGUEES...
\def\ab{\overline{a}}
\def\bb{\overline{b}}
\def\cb{\overline{c}}
\def\db{\overline{d}}
\def\zb{\overline{z}}
\def\ub{\overline{u}}
\def\vb{\overline{v}}
\def\wb{\overline{w}}
\def\fb{\overline{F}}
\def\pb{\overline{p}}
\def\rb{\overline{r}}
\def\xb{\overline{x}}
\def\Hb{\overline{H}}
\def\Ub{\overline{U}}
\def\Wb{\overline{W}}
\def\uh{\hat{u}}

%AUTRES...
\def\qed{\Box}
\def\lvs{\vskip2mm}
\def\no{\mid\!\mid}
\def\sq{\Box_u}
\def\gsq{\Box}

%THEOREMES...
\newtheorem{defi}{Definition}[section]
\newtheorem{prop}{Proposition}[section]
\newtheorem{theo}{Theorem}[section]
\newtheorem{lemm}{Lemma}[section]
\newtheorem{coro}{Corollary}[section]

\noindent {\Large\bf Hopf algebras and subfactors associated to vertex models}

\bigskip\noindent {\large Teodor Banica}

\bigskip\noindent Institut de Math\'ematiques de Luminy,
Marseille, and

\noindent Institut de Math\'ematiques de Jussieu, case 191,
Universit\'e Paris 6, 4 place Jussieu, 75005 Paris, France

\noindent {\em E-mail:} banica@math.jussieu.fr

\bigskip\noindent {\bf Abstract:} If $H$ is a Hopf algebra whose
square of the antipode is the identity, $v\in\l (V)\otimes H$ is a
corepresentation, and $\pi :H\rightarrow\l (W)$ is a representation,
then $u=(id\otimes\pi )v$ satisfies the equation $(t\otimes
id)u^{-1}=((t\otimes id)u)^{-1}$ of the vertex models for
subfactors. A universal construction shows that any solution $u$ of this equation arises in this way. A more elaborate construction shows that there
exists a ``minimal'' triple $(H,v,\pi )$ satisfying $(id\otimes\pi
)v=u$. This paper is devoted to the study of this latter construction
of Hopf algebras. If $u$ is unitary we construct a $\c^*$-norm on $H$
and we find a new description of the standard invariant of the
subfactor associated to $u$. We discuss also the ``twisted'' (i.e. $S^2\neq id$) case.

\section*{Introduction}

Let $V,W$ be finite dimensional linear spaces over a field $k$ and consider the
following condition on an invertible element $u\in\l (V)\otimes \l (W)$
$$(t\otimes id)u^{-1}=((t\otimes id)u)^{-1}\,\,\,\,\, (\star )$$
where $t:\l (V)\rightarrow\l (V^*)$ is the transposition. This
equation appeared in the work of V. Jones, and says that a related
partition function is invariant under the type II Reidemeister
moves. For a first (quite disguised) appearance of $(\star )$ see the formulas (1.3 \& 1.4) in \cite{j}. Any unitary
solution of $(\star )$ may be used for constructing a commuting square
$$\matrix{
\c\otimes\l (W) &\subset & \l (V)\otimes \l (W)\cr 
\cup &\ &\cup\cr 
\c &\subset & u(\l (V)\otimes\c )u^{-1}  \cr
}\,\,\,\,\,\, (\sq )$$
hence two examples of subfactors (see \cite{ksv}; here of course we assume that $k=\c$ and that $V$ and $W$ are Hilbert
spaces). In the general case an extension of this
construction is still available (see \cite{bhj}). For the sake of
completness we recall also that if we take bases such that $V=k^m$
and $W=k^n$ then one can associate a 2d vertex model to any $u\in
M_m(k)\otimes M_n(k)$ in the following way: there are $m$ spins per
vertical edge, $n$ spins per horizontal edge, and $u$ is the matrix of Boltzmann weights. See the book \cite{js} for a global look to these facts.

This paper deals with Hopf algebras and is based on the following
approach to the equation $(\star )$. Consider triples $(H,v,\pi )$
consisting of a Hopf algebra $H$ whose square of the antipode is
the identity, a corepresentation $v\in\l (V)\otimes H$ and a
representation $\pi :H\rightarrow \l (W)$. It is easy to see that
$(id\otimes\pi )v$ satisfies $(\star )$, and a universal
construction shows that any solution of $(\star )$ arises in this
way. More generally, one can describe the class of the elements
of the form $(id\otimes\pi )v$, when no assumption on the square of
the antipode is made - this contains for instance the solutions of some
``twisted'' versions of $(\star )$. Most of the paper is written in
this generality, but for simplicity we restrict now attention to the solutions of $(\star )$.
 
Any solution of $(\star )$ gives rise to a Hopf algebra in the
following way. Let us call models for $u$ the triples $(H,v,\pi )$
such that $(id\otimes\pi )v=u$. We will show that $u$ admits a minimal
model. Here ``minimality'' is by definition a certain universality
property, but we will find several descriptions (including a quite explicit
construction) of the minimal model. It is useful to keep in mind the
following heuristical interpretation: if $(H,v,\pi )$ is a model for
$u$ then $v$ and $\pi$ correspond to representations $G\rightarrow
\g\ll (V)$ and $\widehat{G}\rightarrow\g\ll (W)$, where $G$ is the quantum group represented by $H$ and $\widehat{G}$ is its dual; the minimal model is
then characterised by the fact that these representations are
faithful.

Any unitary solution of $(\star )$ gives rise to a Hopf $\c^*$-algebra
and to a Hopf von Neumann algebra in the following way. If $(H,v,\pi
)$ is the minimal model for $u$ we will construct an involution and a
$\c^*$-norm on $H$ and we will prove that the pair $(\bar{H} ,v)$
satisfies Woronowicz' axioms from \cite{w1},\cite{w2} where $\bar{H}$ is the
completion of $H$. As the square of the antipode is the identity, the
Haar measure is a trace, and by GNS construction one
gets a Kac algebra of compact type in the sense of \cite{es}. It is
useful to keep in mind the following heuristical interpretation:
$\bar{H}=C(G)=\c^*(\Gamma )$, with $G$ a compact quantum group and
$\Gamma =\widehat{G}$ a discrete quantum group. We also show that the
operation $u\mapsto \bar{H}$ produces ``most'' of the commutative and
cocommutative Hopf $\c^*$-algebras - these come from the obvious
solutions $\sum g_i\otimes e_{ii}$ and $\sum e_{ii}\otimes g_i$ of
$(\star )$ - as well as all the finite dimensional ones.

Let us call $L(\sq )$ the standard invariant of the ``vertical'' subfactor associated to the commuting square $\sq$. There is a decription due to V. Jones of this lattice $L(\sq )$ which uses
diagrams (see for instance \cite{js}). We will prove that if $(H,v,\pi
)$ is the minimal model for 
$$u^\prime =u_{12}((id\otimes
t)u^{-1})_{13}\in \l (V)\otimes \l
(W)\otimes \l (W^*)=\l (V)\otimes \l (W\otimes W^*)$$
(which satisfies also $(\star )$) then $L(\sq )$ is equal to the following lattice $L(v)$
$$\matrix{
\c &\subset & End(v) & \subset &
End(v\otimes\hat{v})&
\subset & End(v\otimes\hat{v}\otimes v)&
\subset\
\cdots\cr &\ &\ \cup &\ &\cup &\ &\cup \cr &\ &\ \c & \subset &
End(\hat{v})&
\subset & End(\hat{v}\otimes v)& \subset\
\cdots\cr
 &\ &\ &\ &\cup &\ &\cup \cr &\ &\ \ & \ & \c &
\subset & End(v)&
\subset\ \cdots\cr &\ &\ &\ &\ &\ &\cup \cr
&\ &\ &\ &\ &\ &\cdots &\ \ \ \ \ \cdots\cr
}$$
{\em Comments.} Summing up, our results split the operation $u\mapsto
L(\sq )$ into a composition of four disjoint operations
$$u\mapsto u^\prime\mapsto (H,v,\pi )\mapsto (\bar{H},v)\mapsto L(v)$$
The first operation is very explicit, the second one is the
construction of the minimal model, and in the third one the involution
making $v$ unitary and the maximal $\c^*$-norm are uniquely
determined. About the fourth one, we would like to mention that given
a Popa system $L$, there are at least two reasons for trying to find a
pair $(A,v)$ satisfying Woronowicz' axioms such that $L=L(v)$. First
of all $L(v)$ and its principal graphs have simple interpretations in
terms of representation theory, and in some cases (e.g. when $A$
happens to be commutative or cocommutative) such an equality $L=L(v)$
is very close to the ultimate result in the ``computation'' of $L$. A
second reason is that certain analytical notions like amenability are
supposed to be better understood for Woronowicz algebras - which have a
Haar measure and all the related structures - than for Popa systems or than for subfactors. See the paper \cite{b2} for an introduction to
the lattices of the form $L(v)$, from a point of view close to the one of
\cite{p}.
\lvs
{\em Acknowledgements.} We would like to thank Vaughan Jones for
pointing out that the diagrammatic picture shows that $L(\sq )$
satisfies the axioms for lattices of the form $L(v)$; this led us to
work out the present ``best'' construction $u\mapsto (A,v)$ such that $L(\sq )=L(v)$. We are also grateful to Patrick Polo for useful discussions on Hopf $k$-algebras.

%\vfill\eject 

\section{Construction of the minimal model}

Let $k$ be a field. All the $k$-algebras will have units and the
morphisms between them will be unital. The $k$-coalgebras will
have counits and the morphisms between them will
be counital. Recall that if $C$ is a coalgebra then $C^*$
has a canonical structure of algebra. Conversely, if $A$ is an
algebra then the subspace $A^\circ\subset A^*$ consisting of
linear forms $f$ such that $ker(f)$ contains a finite
codimensional ideal of $A$ has a canonical structure of
coalgebra. Note that if $A$ is finite dimensional then
$A^\circ =A^*$. If $H$ is a Hopf algebra it follows that
$H^\circ$ has a canonical structure of Hopf algebra. See
\cite{abe}.

Let $(H,m,1,\Delta ,\varepsilon ,S)$ be a Hopf $k$-algebra. We
call finite dimensional representations of $H$ the morphisms
of algebras
$\pi :H\rightarrow \l (V)$, where $V$ is a finite
dimensional $k$-linear space. The space of coefficients of
$\pi$ is the following linear subspace of $H^\circ$, which is
easily seen to be a subcoalgebra:
$$\co_\pi =\{ f\pi\mid  f\in \l (V)^*\}\subset H^\circ$$
The dual notion is that of a morphism of coalgebras $\nu
:\l (V)^*\rightarrow H$, but we prefer
to work with the corresponding element in
$\l (V)\otimes H$. That is, we call finite
dimensional corepresentations of $H$ the elements
$v\in \l (V)\otimes H$ satisfying 
$$(id\otimes\Delta
)v=v_{12}v_{13},\,\,\,\, (id\otimes\varepsilon )v=1$$
where $V$ is a finite dimensional $k$-linear space. The space
of coefficients of $v$ is the following linear subspace of
$H$, which is easily seen
to be a subcoalgebra:
$$\co_v =\{ (f\otimes id)v\mid  f\in \l (V)^*\}\subset H$$
The following simple facts will be intensively used without
reference. The element $(id\otimes S)v$ is an inverse for $v$ in $\l
(V)\otimes H$ - this follows by considering $(id\otimes E)v$,
with $E=m(S\otimes id)\Delta =m(id\otimes S)\Delta =\varepsilon
(.)1$. If $t:\l (V)\rightarrow \l (V^*)$ is the transposition, then
$(t\otimes S)v$ is a corepresentation of $H$ - this follows from the
fact that $t$ is unital and antimultiplicative, and $S$ is counital
and anticomultiplicative.

\begin{defi}
Let $V,W$ be finite dimensional $k$-linear spaces and let $u\in
\l(V)\otimes\l (W)$. A model for $u$ is a triple $(H,v,\pi )$
consisting of a Hopf algebra
$H$, a corepresentation $v\in \l (V)\otimes H$ and a
representation
$\pi :H\rightarrow \l (W)$ such that $(id\otimes\pi )v=u$.
\end{defi}

{\em Comments.} The identifications of the form $X\otimes
Y\simeq \l (X^*,Y)$ show that one may use the following related
definition: if $C$ is a coalgebra and $A$ is an algebra, a model for a
linear map $\varphi :C\rightarrow A$ is a factorisation of it as
$$C\mathop{\longrightarrow}^\nu H\mathop{\longrightarrow}^\pi A$$
where $H$ is a Hopf algebra, $\nu$ is a morphism of coalgebras, and
$\pi$ is a morphism of algebras. Another related definition may be
found by extending the following interpretation of $(id\otimes\pi )v$
when $H$ is finite dimensional: the duality $H\otimes H^*\rightarrow
k$ gives by transposition a distinguished element $\xi_H\in H^*\otimes
H$, and if $\rho :H^*\rightarrow\l (V)$ is the representation corresponding to
$v$ then $(id\otimes\pi )v$ is the image of $\xi_H$ by the
representation $\rho\otimes\pi$. For some technical reasons we prefer to use the above Def. 1.1.
\lvs 
If $V$ is a linear space and $S\subset V$ is a subset, the
orthogonal of $S$ is $S^\perp =\{ f\in V^*\mid
f(s)=0,\,\forall\, s\in S\}$. Also if $T\subset V^*$ is a
subset we may define
$T^{\perp
  (V)} =\{ v\in V\mid t(v)=0,\,\forall\, t\in T\}$. Given
$f\in V^*$, by letting the family $\{ f+a^\perp\mid a\in V\}$
be a base for a system of neighborhoods of $f$, $V^*$ becomes a
linear topological space. A linear subspace $T\subset V^*$ is
dense in $V^*$ with respect to this topology if and
only if $T^{\perp (V)}=\{ 0\}$. See \cite{abe}.

It $C\subset H$ is a subcoalgebra of a Hopf algebra, we denote
by
$<C>$ the Hopf subalgebra of $H$ generated by $C$. That is, $<C>$ is
by definition the (unital) subalgebra of $H$ generated by the set
$\cup_{k\geq 0}S^k(C)$.

\begin{defi}
A model $(H,v,\pi )$ for $u\in
\l(V)\otimes\l (W)$ is said to be left-faithful if
$<\co_v>=H$; right-faithful if $<\co_\pi >$
is dense in $H^*$; and bi-faithful if it is both
left- and right-faithful.
\end{defi}

Given a model $(H,v,\pi )$ one can construct a left-faithful
model $(H^\prime ,v,\pi^\prime )$ in the following way:
$H^\prime$ is $<\co_v>$ and $\pi^\prime$ is the restriction to
$H^\prime$ of $\pi$. Due to this simple fact, we will oftenly restrict attention
to left-faithful models. 

We define the morphisms $(H_1,v_1,\pi_1 )\rightarrow
(H_2,v_2,\pi_2)$ of left-faithful models to be the Hopf
algebra morphisms $f:H_1\rightarrow H_2$ such that
$(id\otimes f)v_1=v_2$ and $\pi_1=\pi_2f$. By
left-faithfulness, such a morphism (if it exists) is
surjective, and unique. In particular a morphism from
$(H,v,\pi )$ to itself has to be equal to the identity
morphism. Thus given $u$, the category of left-faithful models
for $u$ has at most one universally repelling object, and at
most one universally attracting object. 

\begin{defi}
The universally repelling (resp. attracting) object in the
category of left-faithful models for $u$ is
called the maximal (resp. minimal) model for $u$.
\end{defi}

In this definition we assume of course that the category is non-empty,
and that the object to be defined exists. The result below clarifies
the situation.

\begin{theo}
Let $V,W$ be finite dimensional $k$-linear spaces and let
$u\in \l(V)\otimes\l (W)$. The following conditions are
equivalent:

(i) there exists a model for $u$.

(ii) there exists a maximal model for $u$.

(iii) there exists a minimal model for $u$.

(iv) $u$ satisfies the following sequence of conditions: $u_0:=u$ is invertible, $u_1:=(t\otimes
id)u_0^{-1}$ is invertible, $u_2:=(t\otimes id)u_1^{-1}$ is
invertible, $u_3:=(t\otimes id)u_2^{-1}$ is invertible, etc..

Moreover, if these conditions are satisfied then the minimal
model for
$u$ may be characterised as the unique bi-faithful model for
$u$.
\end{theo}

{\em Proof of $(i) \Longrightarrow (iv)$.} By recurrence, we
have to prove that if there exists a model for $w\in \l
(T)\otimes \l (W)$ then $w$ is invertible and there exists a
model for $(t\otimes id)w^{-1}$. Let $(H,v,\pi )$ be a model
for $w$. As $(id\otimes S)v$ is an inverse for $v$, we get
that $(id\otimes \pi S)v$ is an inverse for $w$. This implies also that
$(t\otimes id)w^{-1}=(t\otimes \pi S)v$. As $(t\otimes S)v$ is a corepresentation of $H$, this shows that $(H,(t\otimes S)v,\pi )$ is a model for $(t\otimes id)w^{-1}$.
 
{\em Proof of $(iv) \Longrightarrow (ii)$.} We take a basis of $V$
which identifies $V=k^s$. Let $F$ be the free $k$-algebra generated by elements $\{ w^n_{a,b}\}$ with $n\geq 0$ and $a,b\in \{ 1,...,s\}$. Then $F$ has a bialgebra structure with 
$$\Delta (w^n_{a,b})=\sum_{1\leq c\leq
s}w^n_{a,c}\otimes w^n_{c,b}\, ,\,\,
\varepsilon (w^n_{a,b})=\delta_{a,b}$$ 
For every $n$ let
$w^n\in M_{s}(F)$ be the matrix having entries
$(w^n_{a,b})_{1\leq a,b\leq s}$. Let $J\subset F$
be the two-sided ideal generated by the relations coming from
identifying the coefficients in the equalities
$$w^n(w^{n+1})^t=(w^{n+1})^tw^n=1$$
for every $n$, and consider the quotient $H=F/J$. Let $v^n=(id\otimes
p)w^n$, where $p:F\rightarrow H$ is the projection. Then $H$ has the
following universal property $(P)$: given any $k$-algebra $A$ and any
sequence of matrices $V^n\in M_s(A)$ such that
$(V^{n+1})^t=(V^n)^{-1}$ for every $n\geq 0$, there exists a (unique) morphism of algebras $f:H\rightarrow A$ such that $(id\otimes
f)v^n=V^n$ for every $n$.

The following equality holds in $H\otimes H$
$$(v^n_{12}v^n_{13})^{-1}=(v^n_{13})^{-1}(v^n_{12})^{-1}=(v^{n+1}_{13})^t(v^{n+1}_{12})^t=(v^{n+1}_{12}v^{n+1}_{13})^t$$
so by applying $(P)$ with $A=H\otimes H$ and $V^n=v^n_{12}v^n_{13}$ we
get a certain morphism $\Delta_H :H\rightarrow H\otimes H$. Also
by applying $(P)$ with $A=k$ and $V^n=1$ we get a morphism
$\varepsilon_H:H\rightarrow k$, and by applying $(P)$ with $A=H^{op}$
and $V^n=(v^{n+1})^t$ we get a morphism ${\cal S}_H:H\rightarrow
H^{op}$. If $j:H^{op}\rightarrow H$ is the canonical map, it is easy
to see that $(H,m,1,\Delta_H,\varepsilon_H,j{\cal S}_H)$ satisfies the
axioms for a Hopf algebra (by verifying each of them on the generators
$v^n_{ab}$). Once again by $(P)$, the conditions on $u$ in the
statement allow us to define a morphism of algebras $\pi :H\rightarrow
\l (W)$ such that $(id\otimes \pi )v^n=u_n$ for every $n$. This shows
that $(H,v^0,\pi )$ is a model for $u$. 

Let $(K,r,\nu )$ be an arbitrary left-faithful model for $u$. We
define a sequence of corepresentations $r^n\in M_s(K)$ by $r^{2n}=(id\otimes
S^{2n})r$ and $r^{2n+1}=(t\otimes S^{2n+1})r$ for every $n\in\n$,
where $t:M_s(k)\rightarrow M_s(k)$ is the transposition. Then
$(r^n)^{-1}=(id\otimes S)r^n=(r^{n+1})^t$ for every $n$, so the
property $(P)$ gives a morphism of Hopf algebras $f:H\rightarrow K$ such
that $(id\otimes f)v^n=r^n$ for every $n$. This shows that $f$ is a
morphism of left-faithful models $(H,v^0,\pi )\rightarrow (K,r,\nu )$. Thus $(H,v^0,\pi )$ is the maximal model for $u$. 
\lvs

The construction below of the minimal model uses the maximal model and the
following simple fact. Assume that ${\cal C}$ is a category such that
for any object $a$ there exists an object $a_1$ and an arrow
$a\rightarrow a_1$ such that for any object $a_2$ and any arrow
$a\rightarrow a_2$ there exists an arrow $a_2\rightarrow a_1$ making
commutative the triangle. Then if ${\cal C}$ has a universally
repelling object, then it has a universally attracting object. We will
show in the next Lemma that the category of left-faithful models for
$u$ has this property. We begin by explaining what the above-mentioned construction $a\mapsto a_1$ is.

If $(H,v,\pi )$ is a left-faithful model for $u$ we define a left-faithful model $(H_1,v_1,\pi_1)$ and a morphism $p_1:(H,v,\pi )\rightarrow
(H_1,v_1,\pi_1)$ in the following way. Consider the space of
coefficients $\co_\pi\subset H^\circ$. Then $\co_\pi$ is a
subcoalgebra of $H^\circ$ and the orthogonal $\co_\pi^{\perp
(H)}$ is the kernel of $\pi$. As $<\co_\pi >$ is a subalgebra
(resp. subcoalgebra) of $H^\circ$, the orthogonal $<\co_\pi
>^{\perp (H)}$ is a coideal (resp. ideal) of $H$, cf.
\cite{abe}, Th. 2.3.6 (i) (resp. Th. 2.3.2 (ii)). Moreover,
from the invariance of $<\co_\pi >$ under the antipode $S^*$
of $H^\circ$ we get the invariance of
$<\co_\pi >^{\perp (H)}$ under $S$, so the quotient 
$$H_1:=H/<\co_\pi >^{\perp (H)}$$
is a Hopf algebra. As $<\co_\pi >^{\perp (H)}$ is contained in
$\co_\pi ^{\perp (H)}$, which is the kernel of $\pi$, we get a
factorisation $\pi =\pi_1p_1$, where $p_1:H\rightarrow H_1$ is the
projection. If we define $v_1:=(id\otimes p_1)v$, then
$(H_1,v_1,\pi_1)$ is a left-faithful model for $u$, and $p_1$ is a
morphism of left-faithful models.

\begin{lemm}
Let $(H,v,\pi )$ be a left-faithful model and construct
$(H_1,v_1,\pi_1)$ and $p_1$ as above. Let $p_2:(H,v,\pi )\rightarrow
(H_2,v_2,\pi_2)$ be a morphism of left-faithful models for $u$. Then
there exists a morphism of left-faithful models $f:(H_2,v_2,\pi_2)\rightarrow (H_1,v_1,\pi_1)$ such that $p_1=fp_2$.
\end{lemm}

{\em Proof.} Let $J$ be the kernel of $p_2$. As $\pi =\pi_2p_2$, the ideal $J$ is contained in the kernel $\co_\pi^{\perp
  (H)}$ of $\pi$; we want to prove that $J$ is contained in the kernel
$<\co_\pi >^{\perp (H)}$ of $p_1$. By dualising  we want to prove
$$\co_\pi\subset
J^\perp\,\Longrightarrow\, <\co_\pi >\subset J^\perp$$
As $J^\perp$ is stable by $S^*$, it contains the set $\cup_{k\geq
  0}(S^*)^k(\co_\pi )$. Now as $J^\perp$ is a subalgebra of $H^*$ (\cite{abe},
Th. 2.3.6 (ii)), it contains the algebra generated
by $\cup_{k\geq
  0}(S^*)^k(\co_\pi )$, which is $<\co_\pi >$. Thus $J\subset
ker(p_1)$, and we get the desired Hopf algebra morphism from
$H_2$ to $H_1$. $\qed$
\lvs 
{\em Proof of $(ii) \Longrightarrow (iii)$.} If $(H,v,\pi)$ is the
maximal model for $u$, then the Lemma 1.1 shows that the above
construction gives the minimal model.

{\em Proof of the last assertion.} Let $(H,v,\pi )$ be the
minimal model for $u$. By applying the above construction we get a certain
left-faithful model $(H_1,v_1,\pi_1)$. As there
exist morphisms from $(H,v,\pi )$ to $(H_1,v_1,\pi_1)$ in both senses, the unicity of morphisms between
left-faithful models shows that these models are isomorphic. In
particular the kernel of the projection $H\rightarrow H_1$, which is
$<\co_\pi >^{\perp (H)}$ by definition, is zero. Thus $(H,v,\pi )$ is bi-faithful.

Now let $(H^\prime ,v^\prime ,\pi^\prime )$ be an arbitrary bi-faithful model for
$u$. By minimality of $(H,v,\pi )$ we get a Hopf algebra
morphism
$p:H^\prime\rightarrow H$ such that $v=(id\otimes p)v^\prime$ and
$\pi^\prime =\pi p$. The
left-faithfulness shows that $p$ is surjective, and the assertion will
follow from the sequence of inclusions
$$ker(p)\subset Im(p^*)^{\perp (H^\prime )}\subset
<\co_{\pi^\prime}>^{\perp (H^\prime )}\subset \{ 0\}$$
The first inclusion is clear: if $x$ is in $ker(p)$, then
$p^*(f)(x)=fp(x)=0$ for any $f\in H^*$. The second one follows by
dualising $<\co_{\pi^\prime}>\subset Im(p^*)$. The third one is the
definition of the right-faithfulness of
$(H^\prime ,v^\prime ,\pi^\prime )$. $\qed$

\section{Relationship with the vertex models}

The simplest way for finding elements $u$ satisfying the conditions in
the Th. 1.1 is to assume that $u_2=u_0$; in this case the infinity of
conditions (iv) in the Th. 1.1 becomes periodic and true. The
condition $u_2=u$ is nothing but the equation $(\star )$ in the
Introduction. In the sections 2,3,4 we restrict attention to
this case, and we will use the following consequence of the Th. 1.1.

\begin{theo}
Let $V,W$ be finite dimensional $k$-linear spaces. If $u\in
\l (V)\otimes\l (W)$ is invertible and if the equality
$$(t\otimes id)u^{-1}=((t\otimes id)u)^{-1}\,\,\,\,\,\,\,\,
(\star)$$ 
holds in $\l (V^*)\otimes\l (W)$, then there exists a minimal
model for $u$. This is the unique bi-faithful model for $u$.
$\qed$
\end{theo}

The next Proposition shows that in the general case,
imposing this important condition is the same as restricting
attention to the Hopf algebras whose square of the antipode is
the identity.

\begin{prop}
Let $V,W$ be finite dimensional $k$-linear spaces.

(i) Let $H$ be a Hopf algebra whose square of the antipode is the
identity, $v\in \l  (V)\otimes H$ be a corepresentation, and
$\pi :H\rightarrow \l (W)$ be a representation. Then
$u:=(id\otimes \pi ) v$ is invertible and satisfies $(\star )$.

(ii) Assume that $u\in \l (V)\otimes \l (W)$ is invertible and
satisfies $(\star )$. If $(H,v,\pi)$ is the minimal model for
$u$, then the square of the antipode of $H$ is the identity.
\end{prop}

{\em Proof. (i)} We know that $(id\otimes S)v$ is the inverse of
$v$. Also $\hat{v}=(t\otimes S)v$ is a corepresentation, so
$(id\otimes S)\hat{v}$ is an inverse for $\hat{v}$. By combining these
results we get $(t\otimes id)v^{-1}=((t\otimes S^2)v)^{-1}$, and as
$S^2=id$ the assertion follows by applying $id\otimes\pi$ to this equality.

{\em (ii)} This may be proved by using $(t\otimes id)v^{-1}=((t\otimes
S^2)v)^{-1}$ and the bi-faithfulness of the minimal model, but one can
do better. Let ${\cal C}$ be the category of left-faithful models for
$u$ such that the square of the antipode of their subjacent Hopf algebras are
the identities. As the property ``$S^2=id$'' is preserved by
surjective morphisms of Hopf algebras, it is enough to prove that
${\cal C}$ is non-empty. Take a basis which identifies $V=k^s$ and let
$H_s$ be the universal $k$-algebra given with generators
$(w_{ij}^0)_{i,j=1,...s}$ and $(w_{ij}^1)_{i,j=1,...s}$ and with the
relations coming from the equalities
$$w^0(w^1)^t=(w^1)^tw^0=w^1(w^0)^t=(w^0)^tw^1=1$$
The same arguments as in the proof of the Th. 1.1 show that there exists a
Hopf algebra structure on $H_s$ (with the antipode given by
$id\otimes S:w^0\leftrightarrow (w^1)^t$) and a representation $\pi$ such
that $(H_s,w^0,\pi )$ is the universally repelling object in ${\cal
  C}$. The details are left to the reader. $\qed$
\lvs 
We give now some examples of minimal models. Recall that if
$G$ is a group then the group algebra $kG$ is the free
$k$-module over
$G$, with the multiplication induced by the one of $G$. It has
a structure of (cocommutative) Hopf
$k$-algebra with $\Delta (g)=g\otimes g$, $\varepsilon
(g)=\delta_{g,e}1_k$ and $S(g)=g^{-1}$, where $e\in G$ is the
unit. If $k(G)$ is the $k$-algebra of all functions from $G$ to
$k$, then we have identifications $kG^*=k(G)$ (as algebras) and
$k(G)^*=k(G)^\circ =kG$ (as coalgebras). The dual of the Hopf algebra
$kG$ is the (commutative) Hopf algebra $R_k(G):=(kG)^\circ\subset k(G)$ consisting
of representative functions on $G$, i.e. coefficients of finite
dimensional representations of $G$. See \cite{abe}. 

Let $G\subset \g\ll (V)$ be a subgroup, with $V$ finite dimensional. The linear map 
$\nu :kG\rightarrow \l (V)$ which maps $g\in kG$ into
$g\in G\subset \g\ll (V)\subset \l (V)$ is a representation of
$kG$, called the fundamental one. We denote by $k[G]\subset
R_k(G)$ the Hopf subalgebra generated by the coefficients of
the representation $G\rightarrow \l (V)$. The space of these coefficients being the image of the transpose $\nu^*:\l (V)^*\rightarrow kG^*$, we get by restriction a coalgebra morphism from $\l (V)^*$ to $k[G]$, hence a corepresentation $v\in\l (V)\otimes k[G]$, called the fundamental one. These constructions depend of course on the embedding $G\subset \g\ll (V)$.

\begin{prop}
Let $V$ be a finite dimensional space and let $g_1,...,g_n$ be
elements of $\g\ll (V)$. Let $G\subset \g\ll (V)$ be the group
generated by $g_1,...,g_n$.

(i) The element $u=\sum e_{ii}\otimes g_i\in \l (k^n)\otimes
\l (V)$ satisfies $(\star )$. If $w=\sum e_{ii}\otimes
g_i\in \l (k^n)\otimes kG$ and $\nu$ is the fundamental
representation of $kG$, then $(kG,w,\nu )$ is the minimal model
for $u$.

(ii) The element $s=\sum g_i\otimes e_{ii}\in \l (V)\otimes
\l (k^n)$ satisfies $(\star )$. If $v$ is the fundamental
corepresentation of $k[G]$ and $\pi :k[G]\rightarrow \l
(k^n)$ is the linear map $f\mapsto\sum e_{pp}f(g_p)$, then
$(k[G],v,\pi )$ is the minimal model for $s$.
\end{prop}

{\em Proof.} It is clear from definitions that $(k[G],v,\pi )$
and $(kG,w,\nu )$ are left-faithful models for $s$,
respectively $u$.

(i) Let $(H,w_1,\nu_1)$ be the minimal model for $u$ and
$p:kG\rightarrow H$ be the corresponding projection. If $K$ denotes the image of $G\subset (kG)^\times$ by the undelying group morphism $p^\times :(kG)^\times\rightarrow H^\times$, then the
elements of $K$ are group-like elements of $H$, so they are linearly
independent (\cite{abe}, Th. 2.1.2). On the other hand these elements
generate $H$ as a linear space, and it follows that $H$ may be
identified with the group algebra $kK$. Now
the equality $\nu =\nu_1p$ reads $\nu^\times\!\!\mid_G=\nu_1^\times
p^\times\!\!\mid_G$, and as $\nu^\times\!\!\mid_G$ is injective we get
that $p^\times\!\!\mid_G$ is injective. Thus the surjection
$G\rightarrow K$ is an isomorphism, so $p$ is an isomorphism. 

(ii) By transposing the inclusions $k[G]\subset {\cal R}_k(G)\subset
k(G)$ we get surjections $k(G)^*\rightarrow {\cal R}_k(G)^*\rightarrow
k[G]^*$. If $q$ denotes the projection from $kG=k(G)^*$ to $k[G]^*$, then
$\pi^*=qw$. As $<\co_w>=kG$, it follows that
$<\co_\pi >=k[G]^*$. In particular $(k[G],v,\pi )$ is
right-faithful. $\qed$
\lvs
We end with a simple Lemma which gives some more examples.

\begin{lemm}
Assume that $u\in\l (V)\otimes \l (W)$ satisfies $(\star )$.

(i) $\bar{u}=(id\otimes t)u^{-1}$ is equal to $((id\otimes t)u)^{-1}$ and
satisfies $(\star )$. 

(ii) $\hat{u}=(t\otimes id)u^{-1}$ satisfies $(\star )$.

(iii) $u_{12}\bar{u}_{13}u_{14}...\in \l (V)\otimes \l (W\otimes
W^*\otimes W...)$ ($i$ terms) satisfies $(\star )$.

(iv) $u_{1,i+1}\hat{u}_{2,i+1}u_{3,i+1}...\in \l (V\otimes V^*\otimes V...)\otimes\l (W)$ ($i$ terms) satisfies $(\star )$.
\end{lemm}

{\em Proof.} The first equality follows by applying the antimorphism
$t\otimes t$ to the equation $(\star )$. The fact that the four
elements in the statement satisfy $(\star )$ is elementary, but we will give a nice proof which will be used later on. By the Prop. 2.1
(i) it is enough to construct models for them such that the square of
the antipode is the identity. Let $(H,v,\pi )$ be the minimal model
for $u$; by the
Prop. 2.1 (ii) the square of the antipode of $H$ is the identity. If
$\hat{v}=(t\otimes S)v$ and $\hat{\pi}=t\pi S$ then $(H,v,\hat{\pi})$ is a model for $\bar{u}$, $(H,\hat{v},\pi )$ is a model for $\hat{u}$, and 
$$(H,v,(\pi\otimes \hat{\pi}\otimes\pi\otimes
...)\Delta^{(i-1)}),\,\,\, (H,v_{1,i+1}\hat{v}_{2,i+1}v_{3,i+1}...,\pi )$$
($i$ terms) are models for $u_{12}\bar{u}_{13}u_{14}...$
and for $u_{1,i+1}\hat{u}_{2,i+1}u_{3,i+1}...$. $\qed$

\section{Tensor products and intertwiners}

If $({\cal C},\otimes )$ is a $k$-linear tensor category
and $X,\hat{X}$ are objects of ${\cal C}$ we may define the
following lattice of $k$-algebras
$$\matrix{
k &\subset & End(X) & \subset &
End(X\otimes\hat{X})&
\subset & End(X\otimes\hat{X}\otimes X)&
\subset\
\cdots\cr &\ &\ \cup &\ &\cup &\ &\cup \cr &\ &\ k & \subset &
End(\hat{X})&
\subset & End(\hat{X}\otimes X)& \subset\
\cdots\cr &\ &\ &\ &\cup &\ &\cup \cr &\ &\ \ & \ & k &
\subset & End(X)&
\subset\ \cdots\cr &\ &\ &\ &\ &\ &\cup \cr
&\ &\ &\ &\ &\ &\cdots &\ \ \ \ \ \cdots\cr
}$$
where the inclusions are the obvious ones. In our examples the
object $\hat{X}$ will be always the ``dual'' of $X$ in some
canonical sense, so this lattice will be denoted simply
by $L(X)$. See the Prop. 1.1 in \cite{b2} for some precise ``duality
conditions'' to be imposed on $(X,\hat{X})$ as to get Jones
projections, traces, etc..

{\em Example.} Let $H$ be a Hopf $k$-algebra whose square
of the antipode is the identity and consider the tensor category of finite dimensional corepresentations of $H$. That
is, the objects are the finite dimensional corepresentations
of $H$, and if $v\in\l (V)\otimes H$ and $w\in\l (W)\otimes H$
are two corepresentations, then the arrows between them are
the intertwining operators
$$Hom(v,w)=\{ T\in\l (V,W)\mid (T\otimes 1_H)v=w(T\otimes
1_H)\}$$
and their tensor product is $v\otimes w:=v_{13}w_{23}$. For
any corepresentation $v$ we define $\hat{v}$ to be the
contragradient corepresentation $(t\otimes S)v=(t\otimes id)v^{-1}$,
and we use the above notation $L(v)$. 

{\em Example.} Let $A$ be a $k$-algebra. We define a tensor category in the following
way. The objects are elements of $\l (V)\otimes A$, with
$V$ ranging over all finite dimensional linear spaces.
If $v\in\l (V)\otimes A$ and $w\in\l (W)\otimes A$ are two
objects the space of arrows between them is
$$Hom(v,w)=\{ T\in\l (V,W)\mid (T\otimes 1_A)v=w(T\otimes
1_A)\}$$
and their tensor product is $v\otimes w:=v_{13}w_{23}$. If
$u\in\l (V)\otimes A$ is an invertible element satisfying the
condition
$$(t\otimes id)u^{-1}=((t\otimes id)u)^{-1}$$
we define $\hat{u}\in \l (V^*)\otimes A$ to be $(t\otimes
id)u^{-1}$ and we use the notation $L(u)$. 

Note that if $A$ is a Hopf algebra whose square of the antipode is the
identity then the category defined
in the first example is a subcategory of the one in the
second example, and that if $v$ is a corepresentation of
$A$ then the two $\hat{v}$'s (hence the two $L(v)$'s) defined
in these ways are equal. Note also that if $A=\l (W)$ then the
above condition on $u$ is exactly the condition $(\star )$.

\begin{prop}
Let $V,W$ be finite dimensional $k$-linear spaces and assume
that $u\in \l(V)\otimes\l (W)$ satisfies $(\star )$. Then the
element
$$u^\prime=u_{12}((id\otimes t)u^{-1})_{13}\in \l (V)\otimes \l
(W)\otimes \l (W^*)=\l (V)\otimes \l
(W\otimes W^*)$$
satisfies $(\star )$. If $(H,v,\pi )$ is a minimal model
for $u^\prime$ then $L(u^\prime )=L(v)$.
\end{prop}

{\em Note.} The notion of equality for lattices in this statement is
the obvious one: there exists a family of inclusion-preserving
isomorphisms of algebras between the algebras of the first lattice and
the algebras of the second lattice. In fact we will prove a little
more finer statement: the lattices $L(u^\prime )$ and $L(v)$ are equal
as sublattices of $L(V)$, where $L(V)$ is the lattice associated to
the linear space $V$ in the sense of \cite{b2}, i.e. $L(V)$ is the
lattice obtained by applying the above construction with $({\cal C},\otimes )$=
the tensor category of finite dimensional linear spaces, and with
$\hat{V}:=V^*$.
\lvs 
{\em Proof.} The fact that $u^\prime$ satisfies $(\star )$ follows
from the Lemma 2.1. Recall that the argument in there was that if
$(K,w,\nu )$ is the minimal model for $u$, then
$$(K,w,(\nu\otimes t\nu S)\Delta )$$
is a model for $u^\prime$. Let $p:K\rightarrow H$ be the Hopf algebra
surjection constructed by the universal property of the minimal
model, which is such that $\pi p=(\nu\otimes t\nu S)\Delta$. By transposing this equality we get 
$$p^*\pi^*=m_{K^*}(\nu^*\otimes
S^*\nu^*t^*)$$
The image of the right term is invariant under $S^*$, so it follows
that the image of $p^*\pi^*$ is invariant under the antipode of
$K^\circ$. As $p^*$ is injective, we get that $\co_\pi =Im(\pi^*)$ is
invariant under the antipode of $H^\circ$. By right-faithfulness we obtain that $\co_\pi$ generates a dense subalgebra of $H^*$.

We have to prove that $L(u^\prime )=L(v)$. Let $0\leq j\leq i$ and
consider the algebra of $L(v)$ sitting on the $(i,-j)$ position. This
algebra is of the form $End(r)$, with $r$ a certain tensor product
between $v$'s and $\hat{v}$'s. As $u^\prime =(id\otimes\pi )v$, it
follows that the algebra of $L(u^\prime )$ sitting on the $(i,-j)$ position is
$End((id\otimes\pi )r)$. Thus it is enough to prove the
following general result. $\qed$

\begin{lemm}
Let $H$ be a Hopf algebra and $\pi :H\rightarrow\l (W)$ be a finite
dimensional representation such that $\co_\pi$ generates a dense
subalgebra of $H^*$. If $r\in \l (T)\otimes H$ is a finite dimensional
corepresentation then $End((id\otimes\pi ) r )=End(r)$.
\end{lemm}

{\em Proof.} The inclusion $End(r)\subset End((id\otimes\pi ) r)$ is
clear. For the converse let $x\in End((id\otimes\pi ) r )$ be
an arbitrary element and consider the set
$$S_x=\{ f\in H^*\mid [x,(id\otimes f)r]=0\}$$
We have to prove that $x\in End(r)$, which is
equivalent to $S_x=H^*$. First of all, if
$f,g$ are two elements of the algebra $H^*$ then $(id\otimes
fg)r$ is equal to 
$$(id\otimes f\otimes g)(id\otimes\Delta )r=
(id\otimes f\otimes g)(r_{12}r_{13})=((id\otimes
f)r)((id\otimes g)r)$$
and this shows that $S_x$ is a subalgebra of $H^*$. On the
other hand we have $\co_\pi\subset S_x$, so $S_x$ contains the
subalgebra of $H^*$ generated by $\co_\pi$; thus $S_x$ is dense in $H^*$. Let now $e_1,...,e_s$ be a basis
of $T$ and identify $\l (T)$ with $M_s(k)$. By writing $r=\sum
e_{ij}\otimes r_{ij}$ and $x=\sum x_{ij}e_{ij}$ we see that an
element $f\in H^*$ is in $S_x$ if and only if 
$$\sum_k x_{ik}f(r_{kj})=\sum_k f(r_{ik})x_{kj}$$
for any $i,j$. By denoting $a_{ij}=\sum_k x_{ik}r_{kj}-
x_{kj}r_{ik}$ we have shown that
$$S_x=\bigcap_{i,j}a_{ij}^\perp$$
This shows that $S_x$ is an open subspace of $H^*$, so $S_x=H^*$. $\qed$

%\vfill\eject 

\section{Relationship with the commuting squares}

In this section $k$ is a field of characteristic zero. Consider a
Markov commuting square of multimatrix $k$-algebras (see \cite{bhj})
$$\matrix{
A_{01} &\subset & A_{11}\cr 
\cup &\ &\cup\cr 
A_{00} &\subset & A_{10}\cr
}\,\,\,\,\, (\gsq )$$
By performing the basic constructions in all the possible directions we obtain an infinite lattice of multimatrix algebras $(A_{ij})_{i,j\geq 0}$, where the labeling is the obvious one. For every $i\geq 0$ denote by $A_{i\infty}$ the inductive limit $lim_{j\rightarrow\infty}(A_{ij})$. The inclusion $A_{0\infty}\subset A_{1\infty}$ is called the vertical inclusion associated to the initial commuting square. The Jones tower of this inclusion is 
$$A_{0\infty}\subset A_{1\infty}\subset A_{2\infty}\subset A_{3\infty}\subset ...$$
Also for every $0\leq j\leq i$ the canonical inclusion $A_{i0}\subset
A_{i\infty}$ induces an isomorphism $A_{j1}^\prime\cap A_{i0}\simeq
A_{j\infty}^\prime\cap A_{i\infty}$, so the lattice $L(\gsq )$ of higher relative
commutants of the vertical inclusion is
$$\matrix{
A_{01}^\prime\cap A_{00} & \subset & A_{01}^\prime\cap
A_{10} & \subset & A_{01}^\prime\cap A_{20} & \subset
&  A_{01}^\prime\cap A_{30} & \subset\cdots\cr
&\ &\ \cup &\ &\cup &\ &\cup \cr 
&\ & A_{11}^\prime\cap
A_{10} & \subset &
A_{11}^\prime\cap A_{20} &
\subset & A_{11}^\prime\cap A_{30}& \subset\
\cdots\cr &\ &\ &\ &\cup &\ &\cup \cr 
&\ &\ \ & \ & A_{21}^\prime\cap A_{20} &
\subset & A_{21}^\prime\cap A_{30}&
\subset\ \cdots\cr &\ &\ &\ &\ &\ &\cup \cr
&\ &\ &\ &\ &\ &\cdots &\ \ \ \ \ \cdots\cr
}$$
These are analogues of well-known results from the theory of commuting
squares of finite dimensional von Neumann algebras, including the
Compactness Theorem of Ocneanu. See \cite{ghj} for commuting squares,
\cite{oc},\cite{js} for Ocneanu's theorem, and \cite{bhj} for the above analogues.  For the purposes of this paper, one may take the above description of $L(\gsq )$ as a definition for it.

An example of commuting squares is the following one. Let $V,W$ be finite dimensional $k$-linear spaces and let $u\in \l (V)\otimes \l (W)$ be an invertible element. Then the diagram
$$\matrix{
k\otimes \l (W) &\subset & \l (V)\otimes \l (W)\cr 
\cup &\ &\cup\cr 
k &\subset & u(\l (V)\otimes k)u^{-1}  \cr
}\,\,\,\,\, (\sq )$$
is a Markov commuting square if and only if $u$ satisfies the
condition $(\star )$ (see \cite{bhj}). This is an analogue of
the well-known result when $V,W$ are finite dimensional Hilbert spaces
and $u$ is unitary, which will be discussed in detail in the next section.

\begin{theo}
Let $u\in \l (V)\otimes \l (W)$ be an element
satisfying $(\star )$. If $(H,v,\pi )$ is the minimal model
for the element
$$u^\prime=u_{12}((id\otimes t)u^{-1})_{13}\in \l (V)\otimes \l
(W)\otimes \l (W^*)=\l (V)\otimes \l
(W\otimes W^*)$$
then the lattice $L(\sq )$ is equal to the lattice $L(v)$.
\end{theo}

{\em Proof.} This will follow from $L(\sq )=L(u^\prime )$ and from the
Prop. 3.1. Most of the partial results that we need claim no
originality and may be deduced from \cite{ksv},\cite{js} by changing
the ground field, the notations etc.. We prefer to give a
short self-contained proof of $L(\sq )=L(u^\prime )$, using our notations.

{\em Step I.} If $T$ is a finite dimensional linear space we
denote by $e^T$ the Jones projection in $\l (T\otimes T^*)$,
i.e. the element $n^{-1}\sum e_{ij}\otimes e^*_{ij}$
which does not depend on the basis $e_1,...,e_n$ of $T$. For any $i\geq 0$ we define an algebra
$$A^i =\l (V\otimes V^*\otimes V\otimes V^*\otimes ...)$$
($i$ terms). Let us prove that the following diagram is a
sequence of basic basic constructions for commuting
squares (note that the first one is $\sq$)
$$\matrix{
\l (W) &\subset &A^1\otimes \l (W) &\subset &A^{2}\otimes\l (W) &\subset & A^{3}\otimes\l (W) &\subset ...\cr
\cup &\ &\cup &\ &\cup &\ &\cup \cr
k &\subset &ad(U^1)A^{1} &\subset &ad(U^2) A^{2} &\subset &ad(U^3) A^3 &\subset ...
}$$
where $U^i\in A^i\otimes \l (W)=\l (V)\otimes\l (V^*)\otimes
...\otimes \l (W)$ is the element
$$U^i=u_{1,i+1}\hat{u}_{2,i+1}u_{3,i+1}\hat{u}_{4,i+1}...$$
($i$ terms). Here of course the inclusions in the upper line are the
obvious ones, and the Jones projection for the $i$-th inclusion of
this line is by definition $id_{A_{i-2}}\otimes e^T\otimes id_W$,
where $T=V$ if $i$ is odd and $T=V^*$ if $i$ is even. By recurrence, we have to prove that in the diagram below, if the square on the left is a commuting square, then the square on the right is obtained by basic construction from the one on the left.
$$\matrix{
A^i\otimes \l (W) &\subset &A^i\otimes \l (T)\otimes \l (W) &\subset & A^i\otimes \l (T)\otimes \l (T^*)\otimes \l (W)\cr 
\cup &\ &\cup &\ &\cup\cr 
ad(U^i)A^i &\subset & ad(U^{i+1})(A^i\otimes \l (T))&\subset 
& ad(U^{i+2})(A^i\otimes \l (T)\otimes
\l (T^*)) }$$
Here $T=V$ or $V^*$ depending on the parity of $i$. Three of the
algebras in the square on the right are the good ones, and the
remaining one is the subalgebra $X$ of $A^i\otimes \l (T)\otimes \l
(T^*)\otimes\l (W)$ generated by the image of
$ad(U^{i+1})(A^i\otimes\l (T))$ and by the Jones projection $e^T_{23}$. Let $v$ be equal to $u$ if $i$ is even and equal to $\hat{u}=(t\otimes
id)u^{-1}$ otherwise. Then $U^{i+1}=U^i_{12}v_{23}$ and
$U^{i+2}=U^i_{14}v_{24}\hat{v}_{34}$. By the Lemma 2.1 $v$ satisfies
$(\star )$, and it is easy to see that this implies that $ad(v_{24}\hat{v}_{34})e^T_{23}=e^T_{23}$. By applying $ad(U^i_{14})$ we get 
$$ad(U^{i+2})e^T_{23}=ad(U^i_{14}v_{24}\hat{v}_{34})e^T_{23}=e^T_{23}$$
On the other hand the image of the embedding
$$ad(U^i_{13}v_{23})(A^i\otimes
\l (T))\subset A^i\otimes \l (T)\otimes \l (T^*)\otimes\l (W)$$
is $ad(U^i_{14}v_{24})(A^i\otimes
\l (T))$, which is equal to $ad(U^{i+2})(A^i\otimes \l (T))$. By
combining these two results we get that $ad((U^{i+2})^{-1})X$
is the algebra generated by 
$A^i\otimes \l (T)$ and by $e^T_{23}$. But this algebra is the basic construction in the upper line, so it is equal to $A^i\otimes
\l (T)\otimes \l (T^*)$ as desired.

{\em Step II.} Fix $0\leq j\leq i$ and consider the algebra
$A=\l (...\otimes V\otimes V^*\otimes V\otimes ...)$ such that
$A^j\otimes A=A^i$. That is, the product has $i-j$ terms and
begins with $V$ if $j$ is even and with $V^*$ otherwise. Let
also $U\in A\otimes \l (W)$ be such that $U^j\otimes U=U^i$.
That is, 
$$U=v_{1,i-j+1}\hat{v}_{2,i-j+1}v_{3,i-j+1}\hat{v}_{4,i-j+1}...$$
where the tensor product has $j-i$ terms and $v=u$ if $j$ is
even and $v=\hat{u}$ if $j$ is odd. The result in the first step tells
us that the algebra of $L(\sq )$ which sits on the $(i,-j)$ position is 
$$D=(A^j\otimes 1_A\otimes \l (W))^\prime\, \cap\, 
ad(U^i)(A^j\otimes 1_{\l (W)})$$
The commutant of $A^j\otimes 1_A\otimes \l (W)$ is $1_{A^j}\otimes A\otimes 1_{\l (W)}$. By using this and by applying $ad(U^j)^{-1}$ we get an isomorphism
$$D\simeq D^\prime :=(A\otimes 1_{\l (W)})\, \cap\, 
ad(U)(A\otimes 1_{\l (W)})$$
We prove that $x\mapsto x\otimes 1_{\l (W^*)}$ induces an isomorphism 
$$D^\prime\simeq (A\otimes 1_{\l (W)}\otimes 1_{\l (W^*)})\,\cap\, (U_{12}((id\otimes t)U^{-1})_{13})^\prime$$
Indeed, if $T\in A$ then $T\otimes 1_{\l (W)}\otimes 1_{\l (W^*)}$ is in algebra on the right iff
$$(T\otimes 1_{\l (W)}\otimes 1_{\l (W^*)})U_{12}((id\otimes
t)U^{-1})_{13}=
U_{12}((id\otimes
t)U^{-1})_{13}(T\otimes 1_{\l (W)}\otimes 1_{\l (W^*)})$$
By the Lemma 2.1 $(id\otimes t)(U^{-1})=((id\otimes t)U)^{-1}$, so this is equivalent to
$$(U^{-1}(T\otimes 1_{\l (W)})U)_{12}=(((id\otimes
t)U^{-1})(T\otimes 1_{\l (W)})((id\otimes
t)U))_{13}$$
It is easy to see that the equality $X_{12}=Y_{13}$ in $A\otimes
{\l (W)}\otimes {\l (W^*)}$ is equivalent to the existence of $S\in A$ such that $X=S\otimes 1_{\l (W)}$ and $Y=S\otimes 1_{\l (W^*)}$. Thus our condition on $T$ is equivalent to the existence of $S\in A$ such that:
$$(T\otimes 1_{\l (W)})U=U(S\otimes 1_{\l (W)}),\,\, (T\otimes 1_{\l (W^*)})((id\otimes
t)U)=((id\otimes
t)U)(S\otimes 1_{\l (W^*)})$$
The second condition could be obtained by applying $(id\otimes t)$
to the first one. Thus the only condition on $T$ remains the
first one, which is exactly $T\in D^\prime$. 
\lvs 
{\em Step III.} If $S,T$ are linear spaces and $w\in\l (S)\otimes \l (T)$ is invertible we use the notation $\wb =(id\otimes t)w^{-1}$. The Step II gives an isomorphism 
$$D^\prime \simeq End(U_{12}((id\otimes
t)U^{-1})_{13})$$
where the $End$ sign is in the sense of the second section. The element $U_{12}((id\otimes
t)U^{-1})_{13}$ is equal to the product
$$(v_{1,i-j+1}\hat{v}_{2,i-j+1}v_{3,i-j+1}\hat{v}_{4,i-j+1}...)(\overline{v}_{1,i-j+2}\overline{\hat{v}}_{2,i-j+2}\overline{v}_{3,i-j+2}\overline{\hat{v}}_{4,i-j+2}...)$$
By rearranging the terms we get that this is equal to
$$(v_{1,i-j+1}\overline{v}_{1,i-j+2})(\hat{v}_{2,i-j+1}\overline{\hat{v}}_{2,i-j+2})(v_{3,i-j+1}\overline{v}_{3,i-j+2})...$$
If $v^\prime :=v_{12}\overline{v}_{13}$ then $\widehat{v^\prime}=\hat{v}_{12}\overline{\hat{v}}_{13}$, so we get that
$$U_{12}((id\otimes
t)U^{-1})_{13}=v^\prime\otimes \widehat{v^\prime}\otimes v^\prime\otimes \widehat{v^\prime}\otimes ...$$
($i-j$ terms, and where the sign $\otimes$ is in the sense of the
second section). Summing up, we get that $D$ is isomorphic to the
$(i,-j)$ algebra of the lattice $L(u^\prime )$, where $u^\prime
=u_{12}\overline{u}_{13}=u_{12}((id\otimes t)u^{-1})_{13})$. It is
easy to see that these isomorphisms commute with the inclusions, so
$L(\sq )=L(u^\prime )$. The Prop. 3.1 applies and gives the result. $\qed$
\lvs 
{\em Examples.} Let $V$ be a finite dimensional linear space and
let $g_1,...,g_n$ be elements of $\g\ll (V)$.

(i) The element $u=\sum e_{ii}\otimes g_i$ satisfies $(\star )$. Let 
$G$ be the subgroup of $\p\g\ll (V)$ generated by
$\tilde{g}_1,...,\tilde{g}_n$, where $g\mapsto \tilde{g}$ is the
projection $\g\ll (V)\rightarrow \p\g\ll (V)$, and consider the corepresentation 
$$v=\sum e_{ii}\otimes g_i\in M_n(k)\otimes kG$$
Then $L(\sq )=L(v)$. This follows
from the Th. 4.1, from the equality $u_{12}((id\otimes
t)u^{-1})_{13}=\sum e_{ii}\otimes g_i\otimes (g_i^{-1})^t$, from
the Prop. 2.2 (i), and by identifying the subgroup of $\g\ll (V\otimes
V^*)$ generated by the $g_i\otimes
(g_i^{-1})^t$'s with the subgroup of $\p\g\ll (V)$ generated by the $\tilde{g}_i$'s.

(ii) The element $u=\sum g_i\otimes e_{ii}$ satisfies $(\star )$. Let
$G\subset\g\ll (V) $ be the group generated by
$g_2g_1^{-1},...,g_ng_1^{-1}$ and let $v$ be the fundamental
corepresentation of $k[G]$. Then $L(\sq )=L(v)$. This follows from the Th. 4.1, from the equality
$u_{12}((id\otimes t)u^{-1})_{13}=\sum g_ig_j^{-1}\otimes
e_{ii}\otimes e_{jj}$, and from the Prop. 2.2 (ii).

These results were obtained in \cite{bhj}; the reader may find in
there the interpretation of the lattice $L(v)$ and of its principal
graphs for these special kinds of corepresentations $v$.
%\vfill\eject 

\section{Hopf $\c^*$-algebras and (twisted) biunitaries}

In this section $k=\c$. Recall that an involution of a $\c$-algebra
$A$ is a unital antilinear antimultiplicative map $*:A\rightarrow A$
such that $*^2=id$. A $\c^*$-norm on $A$ is a norm making
$A$ into a normed algebra, and such that $\no a^*a\no =\no
a\no^2$ for any $a\in A$. The completion of $A$ with respect to the
biggest $\c^*$-norm (if such a norm exists) is a $\c^*$-algebra called
envelloping $\c^*$-algebra and denoted here $\bar{A}$. An involution
of a Hopf $\c$-algebra $H$ is an involution $*$ of the subjacent
algebra such that for any $a\in H$ the following formulas hold
$$\Delta (a^*)=\Delta (a)^*,\,\varepsilon (a^*)=\overline{\varepsilon
  (a)},\, S(S(a^*)^*)=a$$
If $V$ is a Hilbert space, a corepresentation $v\in\l (V)\otimes H$
  which is a unitary element of the $*$-algebra $\l (V)\otimes H$ is
  said to be unitary corepresentation.

In \cite{w1}, \cite{w2} it was developed a theory for pairs $(A,v)$
consisting of a unital $\c^*$-algebra and a unitary $v\in\l (V)\otimes
A$, where $V$ is some finite dimensional Hilbert space, subject to the
following conditions, to be reffered as ``Woronowicz' axioms''.

\noindent - the coefficients of $v$ generate $A$ as a
dense $*$-subalgebra (called $A_s$).

\noindent - there exists a $\c^*$-morphism $\delta :A\rightarrow
A\otimes_{min}A$ such that $(id\otimes\delta )v=v_{12}v_{13}$.

\noindent - there exists a linear antimultiplicative map $\kappa
:A_s\rightarrow A_s$ such that $(id\otimes\kappa )v=v^{-1}$ and such
that $\kappa (\kappa (a^*)^*)=a$ for
any $a\in A_s$.

These conditions imply that $A_s$ has a canonical structure of Hopf
$*$-algebra, and that $v$ is a unitary corepresentation (see
\cite{w1}). While we are interested only in such pairs $(A,v)$ we recall that
with a suitable choice of arrows, the category ${\cal C}$ of inductive
limits of such objects is called the category of ``Woronowicz
algebras'', or ``unital Hopf $\c^*$-algebras''. By reversing the
arrows of ${\cal C}$ we get the category of ``compact quantum groups''
and by reversing them once again (sic!) we get the category of
``discrete quantum groups''.

\begin{defi}
A twisted biunitary is a unitary $u\in \l (V)\otimes \l (W)$, where
$V$ and $W$ are finite dimensional Hilbert spaces, which satisfies the
equation
$$((Q^t)^{-2}\otimes id)((t\otimes id)u^{-1})((Q^t)^2\otimes
id)=((t\otimes id)u)^{-1}\,\,\,\, (\star_Q)$$
for some positive operator $Q\in\l (V)$.
\end{defi}

For $Q=id$ the condition $(\star_Q)$ is exactly $(\star )$, and is
equivalent to the fact that $\sq$ is a commuting square in the sense
of subfactor theory; in this case $u$ is said to be a biunitary (see \cite{ksv}). By taking suitable bases the twisted biunitarity condition has the following equivalent formulation. A unitary $u=(u^{ax}_{by})\in M_m(\c )\otimes M_n(\c )$ is said to be a twisted biunitary if there exist positive real numbers $q_1,...,q_m$ such that
$$\sum_{b,x}
q_b\overline{u}^{ax}_{by}u^{bz}_{cx}=q_a\delta_{a,c}\delta_{y,z},\,\,
\forall\, a,c,y,z$$

\begin{prop}
If $(A,v)$ satisfies Woronowicz' axioms and $\pi :A\rightarrow\l (W)$
is a $*$-representation then $(id\otimes\pi )v$ is a twisted
biunitary. Any twisted biunitary arises in this way.
\end{prop}

We will need the following easy lemma, which gives in particular some
more equivalent formulations of the twisted biunitarity condition.

\begin{lemm}
If $V$ is a finite dimensional Hilbert space, $Q\in\l (V)$ is a
positive operator, $A$ is a $\c^*$-algebra, and $u\in\l (V)\otimes A$
is a unitary then the following conditions are equivalent:

(i) $(Q^t\otimes id)((t\otimes id)u^{-1})((Q^t)^{-1}\otimes
id)$ is unitary.

(ii) $((Q^t)^{-1}\otimes id)((t\otimes id)u)(Q^t\otimes
id)$ is unitary.

(iii) $((Q^t)^{-2}\otimes id)((t\otimes id)u^{-1})((Q^t)^2\otimes
id)=((t\otimes id)u)^{-1}$. $\qed$
\end{lemm}

{\em Proof of the Prop. 5.1.} Consider the family of characters
$(f_z)_{z\in\c}:A_s\rightarrow\c$ introduced in \cite{w1} and let
$Q=(id\otimes f_{\frac{1}{2}})v$. Then $Q>0$ and 
$$(Q^t\otimes id)((t\otimes id)v^{-1})((Q^t)^{-1}\otimes id)=(t\otimes
j)v$$
where $j:A_s\rightarrow A_s$ is the linear map $x\mapsto
f_{\frac{1}{2}}*\kappa (x)*f_{-\frac{1}{2}}$ and $*$ is the
convolution over $A_s$ (cf. the formulas in the $5^{th}$ section of
\cite{w1}). These formulas show also that $j$ is an antimorphism of
$*$-algebras, so $t\otimes j$ is also an antimorphism of $*$-algebras,
so it maps unitaries to unitaries. By applying $id\otimes\pi$ to $(t\otimes j)v$ we get the result. See the proof of the Lemma 1.5 in
\cite{b2} for more details.

Conversely, let $u\in\l (V)\otimes\l (W)$ be a twisted
biunitary. Choose an orthonormal basis in $V$ consisting of
eigenvectors of $Q$, so that $V=\c^n$ and $Q\in M_n(\c )$ is diagonal and positive. Consider
the universal $\c^*$-algebra $A_u(Q)$ generated by the coefficients of
a unitary $n\times n$ matrix $w$ such that $Q\overline{w}Q^{-1}$ (or,
equivalently, $Q^{-1}w^tQ$) is also unitary. Then $(A_u(Q),w)$
satisfies Woronowicz' axioms (see \cite{b1}) and as $u$ and
$Q^{-1}u^tQ$ are unitaries (cf. the Lemma 5.1) we get a $\c^*$-morphism $f :A_u(Q)\rightarrow\l (W)$ such that
$(id\otimes f )w=u$. $\qed$
\lvs
With the above notations $(A_u(F)_s,w,f)$ is a model for $u$. Thus
the Th. 1.1 applies to any twisted biunitary. In fact the conditions
(iv) of the Th. 1.1 are also easy to verify - one gets by recurrence
that for every $n\geq 0$:
$$u_{2n}=(Q^{2n}\otimes id)u(Q^{-2n}\otimes id),\,\,
u_{2n+1}=((Q^{-2n})^t\otimes id)u_1((Q^{2n})^t\otimes id)$$

\begin{theo}
Let $u\in \l (V)\otimes \l (W)$ be a twisted biunitary. If
$(H,v,\pi )$ is the minimal model for $u$ then
there exists an involution on $H$ such that:

(i) $v$ is a unitary representation and $\pi$ is a
$*$-representation.

(ii) $H$ has a (biggest) $\c^*$-norm.

(iii) $(\bar{H},v)$ satisfies Woronowicz' axioms.

\noindent Moreover, this involution of $H$ is the unique one such that
$v$ is unitary.
\end{theo}

{\em Proof.} Choose a basis such that $V=\c^n$ and $Q\in M_n(\c )$ is
positive and diagonal. Consider the model $(A_u(Q)_s,w,f)$
constructed in the proof of the Prop. 5.1, and let $q:A_u(Q)_s\rightarrow
H$ be the corresponding Hopf algebra morphism. The following equality
holds in $M_n(\c )\otimes A_u(Q)$
$$(id\otimes S^2)w=(t\otimes id)((t\otimes id)w^{-1})^{-1}=(Q^2\otimes id)w(Q^{-2} \otimes id)$$
By applying $id\otimes q$ to this equality we get
$$(id\otimes  S^2)v=(Q^2\otimes id)v(Q^{-2} \otimes id)$$
In particular $\co_v$ is stable under $S^2$, and by left-faithfulness
we obtain that $H$ is generated as an algebra by $\co_v$ and $S(\co_v )$.

Let us prove firstly the unicity part. Such an involution $*$ has to satisfy
$v^*=v^{-1}=(id\otimes S)v$, so it is uniquely determined on
$\co_v$. As $*^2=id$ the restriction of $*$ to $S(\co_v )$ is the
inverse of the restriction of $*$ to $\co_v$, so it is also uniquely
determined. As $\co_v$ and $S(\co_v )$ generate $H$ as an algebra, $*$
extends uniquely by antimultiplicativity.

We denote as usual by $*$ the involutions of $\l (V)$ and $\l (W)$. Let $K$ be
complex conjugate of $H$ and denote by $j:H\rightarrow K$ the canonical
antilinear isomorphism. Then $*\pi j^{-1}$ is a representation
of $K$ and $(*\otimes jS)v$ is a corepresentation
of $K$. As $u=(id\otimes\pi )v$ and $v^{-1}=(id\otimes
S)v$ we get that 
$$u=(u^{-1})^*=(*\otimes *)(id\otimes\pi S)v=(id\otimes *\pi
j^{-1})(*\otimes jS)v$$
so $(K ,(*\otimes jS)v, *\pi j^{-1})$ is a model for $u$. We have
$S^2(\co_v)=\co_v$, so this model is left-faithful and by the
universality property of the minimal model we get a Hopf
algebra morphism $p:K\rightarrow H$ such that $\pi p=*\pi j^{-1}$ and
$(*\otimes pjS)v=v$. With $r=pj$ these formulas are
$$\pi r=*\pi\, ,\,\,\, (*\otimes rS)v=v\,\,\,\,\,\, (\dag )$$
We prove that $r$ is an involution of $H$. The facts that $r$ is
unital, antilinear, antimultiplicative, anticounital and
comultiplicative follow from the corresponding properties of $p$ and
$j$. From $(id\otimes S)v=v^{-1}$ and $(\dag )$ we get that $(*\otimes r)v^{-1}=v$. As
$*\otimes r$ is an antiautomorphism of the algebra $M_n(\c )\otimes H$, this shows also that $(*\otimes r)v=v^{-1}$. By
applying $*\otimes r$ to these two formulas we get that $(id\otimes
r^2)v=v$ and $(id\otimes
r^2)v^{-1}=v^{-1}$. Thus $r^2=id$ on both $\co_v$ and $S(\co_v )$, and
as these spaces generate $H$ as an algebra we get that $r^2=id$. The
proof of $(rS)^2=id$ is similar: we have $(*\otimes rS)v=v$, and as
$*\otimes rS$ is an antiautomorphism we get also that $(*\otimes
rS)v^{-1}=v^{-1}$. By applying $*\otimes rS$ we get that $(rS)^2=id$
on both $\co_v$ and $S(\co_v )$, and this implies that $(rS)^2=id$.

Thus $r$ is an involution of $H$; from now on we denote it by $*$. The
point (i) is clear from the above formulas $(\dag )$. As $v$ is unitary, the norms of its coefficients are less
than one for every $\c^*$-seminorm on $H$. By left-faithfulness
these coefficients generate $H$ as a $*$-algebra, so there exists a maximal
$\c^*$-seminorm $\no\ \no$ on $H$. The point (ii) is equivalent to the
fact that this is a norm. If $A$ denotes the completion of the separation
of $H$ by $\no\ \no$, this is the same as proving that the
canonical map $i:H\rightarrow A$ is injective. 

Consider the $*$-morphism $(i\otimes i)\Delta :H\rightarrow
A\otimes_{min}A$. It has values in a $\c^*$-algebra, so it extends to a $\c^*$-morphism $\delta :A\rightarrow
A\otimes_{min}A$ which satisfies $(id\otimes\delta )V=V_{12}V_{13}$, where $V=(id\otimes
i)v$. On the other hand from the fact that $w^t$ is invertible in
$M_n(A_u(Q))$ we get that $v^t$ is invertible in $M_n(H)$, so $V^t$ is
invertible in $M_n(A)$. Summing up, the pair $(A,V)$
satisfies the first two axioms of Woronowicz and is such that $V^t$ is invertible; by \cite{w3} we get that $(A,V)$
satisfies all Woronowicz' axioms. In particular we get a Hopf algebra
structure on $i(H)=A_s$ such that $V$ is a corepresentation of it.

Consider the $*$-morphism $\pi :H\rightarrow \l (W)$. It has
values in a $\c^*$-algebra, so it extends to a $\c^*$-morphism
$\bar{\pi}:A\rightarrow\l (W)$. Thus $(A_s,V,\bar{\pi}_{\mid A_s})$
is a model for $u$. By the universality property of $(H,v,\pi )$ we
get a section for $i$, so $i$ is injective. This finishes the proof of (ii) and (iii). $\qed$
\lvs
{\em Example.} Any finite dimensional Hopf $\c^*$-algebra is subjacent
to a minimal model for a biunitary. Indeed, if $A$ is finite
dimensional, then the square of its antipode is the identity (see
\cite{w1}), and $(A,v,\pi )$ is a bi-faithful model for
$u=(id\otimes\pi )v$, where $\pi$ is the regular representation, and $v$ is the coregular corepresentation.
\lvs 
{\em Example.} Let $V$ be a finite dimensional Hilbert space
and let $g_1,...,g_n$ be elements of $\u (V)$. Let $G\subset\u
(V)$ be the group generated by $g_1,...,g_n$.

(i) The element $u=\sum e_{ii}\otimes g_i$ is biunitary. Consider the
minimal model $(kG,w,\nu )$ for $u$ (cf. Prop 2.2 (i)). Then the involution constructed in the Th. 5.1 is given by $\sum
c_ig_i\mapsto\sum\bar{c}_ig_i$, and the $\c^*$-algebra $\overline{kG}$ is the
(full) group $\c^*$-algebra $\c^*(G)$.

(ii) The element $u=\sum g_i\otimes e_{ii}$ is a biunitary. Consider
the minimal model $(k[G],v,\pi )$ for $u$ (cf. Prop 2.2 (ii)). Then
the involution constructed in the Th. 5.1 is given by $f\mapsto
(g\mapsto \bar{f}(g))$, and the classical Peter-Weyl theory shows that
the $\c^*$-algebra $\overline{k[G]}$ is the algebra $C(\bar{G})$ of
continuous functions on the closure $\bar{G}\subset\u (V)$ of
$G$.

By combining these results with the Th. 4.1 one gets the descriptions
from  \cite{ksv} of the standard invariants and of the principal
graphs of the corresponding subfactors.

\end{document}